\documentclass[10pt]{amsart}
\allowdisplaybreaks[4]
\usepackage{amssymb,color}
\usepackage[left=34mm,right=30mm,top=40mm,bottom=32mm]{geometry}
\usepackage{graphicx}
\usepackage{float}
\usepackage[colorlinks=true, citecolor=blue, linkcolor=blue]{hyperref}
\usepackage{amsthm}

\definecolor{c20}{rgb}{0.,0.7,0.}
\definecolor{c30}{rgb}{0.,0.,1.}
\definecolor{c40}{rgb}{1,0.1,0.7}
\definecolor{c50}{rgb}{1,0,0}
\definecolor{c60}{rgb}{0,0.9,0.1}

\def\xH#1{\textcolor{c20}{#1}}
\def\xH#1{#1}
\def\xx#1{\textcolor{c60}{#1}}
\def\xx#1{#1}

\def\cJ#1{\textcolor{c50}{#1}}
\def\cJ#1{#1}

\def\cL#1{\textcolor{c50}{#1}}
\def\cL#1{#1}
\def\cc#1{\textcolor{c50}{#1}}
\def\cc#1{#1}

\newcommand{\abs}[1]{\left\lvert #1 \right\rvert}

\newcommand{\E}[1]{\mathbb{E}\left \{#1\right\}}

\newcommand{\pk}[1]{\mathbb{P} \left \{#1 \right\} }

\newcommand{\R}{\mathbb{R}}

\newcommand{\BQN}{\begin{eqnarray}}
\newcommand{\EQN}{\end{eqnarray}}
\newcommand{\BQNY}{\begin{eqnarray*}}
\newcommand{\EQNY}{\end{eqnarray*}}

\newcommand{\BS}{\begin{sat}}
\newcommand{\ES}{\end{sat}}
\newcommand{\BT}{\begin{theo}}
\newcommand{\ET}{\end{theo}}
\newcommand{\BK}{\begin{korr}}
\newcommand{\EK}{\end{korr}}

\newcommand{\BD}{\begin{de}}
\newcommand{\ED}{\end{de}}

\newcommand{\BIT}{\begin{itemize}}
\newcommand{\EIT}{\end{itemize}}
\newcommand{\BDI}{\begin{description}}
\newcommand{\EDI}{\end{description}}

\newcommand{\BRM}{\begin{remarks}}
\newcommand{\ERM}{\end{remarks}}

\newcommand{\BEL}{\begin{lem}}
\newcommand{\EEL}{\end{lem}}

\newtheorem{theo}{Theorem}[section]
\newtheorem{sat}[theo]{Proposition}
\newtheorem{de}[theo]{Definition}
\newtheorem{lem}[theo]{Lemma}

\newtheorem{korr}[theo]{Corollary}
\newtheorem{remark}[theo]{Remark}
\newtheorem{remarks}[theo]{Remarks}

\newcommand{\nelem}[1]{{Lemma \ref{#1}}}

\newcommand{\netheo}[1]{{Theorem \ref{#1}}}

\newcommand{\prooftheo}[1]{ \textbf{Proof of Theorem} \ref{#1} }

\newcommand{\prooflem}[1]{\textbf{Proof of Lemma} \ref{#1}}

\newcommand{\COM}[1]{}

\newcommand{\QED}{\hfill $\Box$ \\}

\def\Ga{\gamma}
\def\vp{\varepsilon}
\def\vn{\varepsilon}

\def\rw{\rightarrow}

\def\IF{\infty}

\def\LT{\left}
\def\RT{\right}

\begin{document}

\title{Ruin problem for Brownian motion risk model with interest rate and tax payment}

\thispagestyle{empty}

{\author{Long Bai}
\address{Long Bai, Department of Actuarial Science, University of Lausanne, UNIL-Dorigny 1015 Lausanne, Switzerland}
\email{long.bai@unil.ch}

\author{Peng Liu}
\address{Peng Liu, Department of Actuarial Science, University of Lausanne,\\
	UNIL-Dorigny, 1015 Lausanne, Switzerland
}
\email{peng.liu@unil.ch}}

\bigskip

\date{\today}
 \maketitle

\bigskip
{\bf Abstract:} Let $\{B(t), t\ge 0\}$ be a  Brownian motion. Consider the Brownian motion risk model with interest rate collection and tax payment defined by
\BQN\label{Rudef}
\widetilde{U}_\Ga^\delta(t)=\widetilde{X}^\delta(t)-\Ga\sup_{s\in[0,t]}
\left(\widetilde{X}^\delta(s)e^{\delta(t-s)}-ue^{\delta(t-s)}\right),t\ge0,\EQN
with $$\widetilde{X}^\delta(t)=ue^{\delta t}+c \int_0^t e^{\delta (t-v)}dv-\sigma \int_0^t e^{\delta(t-v)}dB(v),$$
where $c>0,\gamma \in [0,1)$ and $\delta \in \R$ are three given constants. When $\delta=0$ and $\gamma \in (0,1)$ this is the risk model introduced from Albrecher and Hipp in \cite{AH2007} where the ruin probability in the infinite time horizon has been explicitly calculated. In the presence of interest rate $\delta\neq0$, the calculation of ruin probability for this risk process for both finite and infinite time horizon seems impossible. In the paper, based on asymptotic theory we propose an approximation for ruin probability and ruin time when the initial capital $u$ tends to infinity. Our results are of interest given the fact that this can be used as benchmark model in various calculations.\\
{\bf Key Words:} Brownian motion; force of interest; tax payment; ruin probability; ruin time.\\
{\bf AMS Classification:} Primary 60G15; secondary 60G70

%%%%%%%%%%%%%%%%%%%%%%%%%%%%%%%%%%%5555555
\section{Introduction}
The risk reserve process of an insurance company without interest can be modelled by a stochastic process $\{ X(t),t\ge0\}$ given as
$$X(t)=u+ct-\sigma B(t),\ t\ge0,$$
see \cite{MR1458613,rolski2009stochastic}, where $u\ge0$ is the initial reserve, $c>0$ is the rate of premium received by the insurance company, and $\sigma B(t)$ is frequently referred as the loss rate of the insurance company.
If we add the effect of tax into the model, the new {\it claim surplus process} is
$$U_\gamma(t)=X(t)-\Ga\sup_{s\in[0,t]}\left(X(s)-u\right),\ t\ge0,$$
where $\gamma\in[0,1]$ is referred to the rate of tax. One of the most important characteristics in risk theory is the ruin probability and \cite{AH2007} shows that
$$\pk{\inf_{t\in[0,T]}U_\gamma(t)<0}=
1-\LT(1-e^{-\frac{cu}{\sigma^2}}\RT)^{\frac{1}{1-\gamma}}.$$
%This exact asymptotic has already dealt with in \cite{HA2013}.\\
Due to the nature of the financial market, we shall consider a more general surplus process including interest rate, see \cite{rolski2009stochastic}, called a risk reserve process with constant force of interest and tax,  i.e., $\widetilde{U}_\Ga^\delta(t),\ t\geq 0$, in \eqref{Rudef}.
%In this paper we consider the effect of interest which plays a pivotal role in the real world, then the model will be
%$$\widetilde{U}_\Ga^\delta(t)=\widetilde{X}^\delta(t)-\Ga\sup_{s\in[0,t]}
%\left(\widetilde{X}^\delta(s)e^{\delta(t-s)}-ue^{\delta(t-s)}\right),\ t\ge0,$$
%where $$\widetilde{X}^\delta(t)=ue^{\delta t}+c \int_0^t e^{\delta (t-v)}dv-\sigma \int_0^t e^{\delta(t-v)}dB(v).$$
For $T\in(0,\IF]$, this contribution is concerned with the exact extreme of $$\psi^\delta_{\Ga,T}(u):=\pk{\inf_{t\in[0,T]}\widetilde{U}_\Ga^\delta(t)<0},$$
as $u\to\IF$. See \cite{rolski2009stochastic, DHJ13a,HX2007} for more studies on risk models with force of interest. Figure \ref{fig1} in Appendix depicts the ruin scenario.
 \\
When $\gamma=0$ and $\delta>0$, i.e. the risk reserve process with positive constant force of interest but without tax, $\psi^\delta_{0,T}(u)$ with $T\in(0,\IF)$ is investigated in \cite{ParisianBrownianfinite2017} and $\psi^\delta_{0,\IF}(u)$ is derived in \cite{Threshold2016,ParisianInfinite2018}.\\
Complementary, we investigate the asymptotic properties of the first passage time (ruin time) of $\widetilde{U}_\Ga^\delta(t)$ on the time interval $[0,T]$, given the ruin has ever  happened during  $[0,T]$. For any $u\geq 0$, and any $T\in(0,\IF]$, define the ruin time of the risk process $\widetilde{U}_\Ga^\delta(t)$ by
\BQN\label{tu}
\tau(u)=\inf\{t\geq 0:\widetilde{U}_\Ga^\delta(t)<0\}.
\EQN
We are interested in the approximate distribution of $\tau(u)|\tau(u)\leq T$, as $u\rw\IF$.\\
Brief organization of the rest of the paper: In Section 2 we first present our main results on the asymptotics of $\psi^\delta_{\Ga,T}(u)$ as $u\rw\IF$ and then
we display the approximation of the ruin time. All the  proofs are relegated to Section 3.

\section{Main Results}
Note that
\BQNY
\psi^\delta_{\Ga,T}(u)=\pk{\inf_{t\in[0,T]}U_\Ga^\delta(t)<0}
=\pk{\sup_{t\in[0,T]}(u-U_\Ga^\delta(t))>u},
\EQNY
where
\BQNY
U_\Ga^\delta(t)=\widetilde{U}_\Ga^\delta(t)e^{-\delta t},t\ge0.
\EQNY
Thus in the analysis of our main results, we consider $\pk{\sup_{t\in[0,T]}(u-U_\Ga^\delta(t))>u}$.\\
In the following theorem, $ \Psi(u), u\in\R$ denotes the survival function
of the standard normal distribution $N(0,1)$. Throughout this paper we write $ f(u)=h(u)(1+o(1))$ or $f(u) \sim h(u)$ if $ \lim_{u \to \infty} \frac{f(u)}{h(u)} = 1 $
and $ f(u) = o(h(u)) $ if $ \lim_{u \to \infty} \frac{f(u)}{h(u)} = 0 $.\\
We prefer to state our new results first, i.e., $\psi_{\gamma,T}^\delta(u)$ as $u\rightarrow\IF$. Cases $T=\IF$
and $T\in(0,\IF)$ are very different and will therefore be dealt with separately. We shall analyse first the case $T\in(0,\IF)$.
\BT\label{TheomOT}
 We have for $\gamma\in[0,1), \delta\in(-\IF,0)\cup(0,\IF)$ and $T\in(0,\IF)$
\BQN\label{result11}
\psi_{\gamma,T}^\delta(u) \sim \frac{2(1+e^{-2\delta T})}{1-\gamma+e^{-2\delta T}}\Psi\left(\frac{u+\frac{c}{\delta}(1-e^{-\delta T})}{a}\right),\ u\rw\IF,
\EQN
where $a^2=\frac{\sigma^2}{2\delta}(1-e^{-2\delta T}).$
\ET
\begin{remarks}
i) When $\gamma=0$, the result of \textbf{\netheo{TheomOT}} reduce to asymptotic ruin probability of the risk model without tax, i.e.
\BQNY
\psi_{\gamma,T}^\delta(u)\sim 2\Psi\left(\frac{u+\frac{c}{\delta}(1-e^{-\delta T})}{a}\right),
\EQNY
which corresponds to the result in \cite{ParisianBrownianfinite2017}.\\
ii) In \eqref{result11}, note that $\gamma$ is just related to the denominator of the constant part. The asymptotic result is increase as $\gamma$.  In fact, when $\gamma$ is bigger, it means that a company need to pay more tax before the ruin happens, thus the ruin probability should be increasing. Table \ref{tab1} is the simulated asymptotic results of $\psi_{\gamma,T}^\delta(u) $ in \netheo{TheomOT}, which also shows the increasing about $\gamma$.
\end{remarks}
\begin{table}[H]
\centering
\caption{The simulated asymptotic results of $\psi_{\gamma,T}^\delta$}\label{tab1}
\begin{tabular}{p{1cm}p{1cm}p{1cm}p{1.5cm}p{1cm}p{1cm}p{4cm}}
\hline
$u$   & $c$   & $\sigma$ & $\delta$ & $T$  & $\gamma$ & asymptotic results \\
\hline
5   & 0.1 & 1      & 0.05   & 20   &0.1    & 0.0363 \\
\hline
5   & 0.1 & 1      & 0.05   & 20  &0.2    & 0.0402 \\
\hline
5   & 0.1 & 1      &-0.05   & 20  &0.2    & 0.0455 \\
\hline
5   & 0.1 & 1      & 0.07   & 20   &0.1    & 0.0210 \\
\hline
5   & 0.1 & 1      & 0.07   & 30   &0.2    & 0.0229 \\
\hline
5   & 0.1 & 1      & -0.07   & 30   &0.2    & 0.0349 \\
\hline
5   & 0.1 & 1      & 0.1   & 20    &0.1   & 0.0090 \\
\hline
5   & 0.1 & 1      & 0.1   & 30    &0.2   & 0.0096 \\
\hline
5   & 0.1 & 1      & -0.1   & 30    &0.2   & 0.0136 \\
\hline
4   & 0.1 & 1      & 0.1   & 20    &0.1   & 0.0312 \\
\hline
4   & 0.1 & 1      & 0.1   & 30    &0.2   & 0.0333 \\
\hline
4   & 0.1 & 1      & -0.1   & 30    &0.2   & 0.0453 \\
\hline
\end{tabular}
%\\*we always set $u=5, c=0.1,\sigma=1, \lambda=600$.
\end{table}

\BT\label{TheomIF} We have for $\gamma\in[0,1), \delta >0$ and $T=\IF$
\BQNY
\psi_{\gamma,\IF}^\delta (u)\sim
\frac{1}{1-\gamma}\widehat{\mathcal{P}}^{\frac{c^2}{\sigma^2\delta}}[0,\IF)
\Psi\LT(\frac{\sqrt{2}}{\sigma}\sqrt{\delta u^2+2cu}\RT),\ u\rw\IF,
\EQNY
where for $-\IF\leq S_1<S_2\leq \IF$
$$
\widehat{\mathcal{P}}^{\frac{c^2}{\sigma^2\delta}}[S_1,S_2]=
\E{\sup_{t\in[S_1,S_2]}e^{\sqrt{2}B(t)-t-\LT(\sqrt{t}-\frac{c}{\sigma\sqrt{\delta}}\RT)^2
}}\in(0,\IF).
$$
\ET
\begin{remarks}
i) We have that
\begin{align*}
\frac{1}{1-\gamma}\widehat{\mathcal{P}}^{\frac{c^2}{\sigma^2\delta}}[0,\IF)
\Psi\LT(\frac{\sqrt{2}}{\sigma}\sqrt{\delta u^2+2cu}\RT)&\sim
\frac{1}{1-\gamma}\widetilde{\mathcal{P}}^{\frac{c^2}{\sigma^2\delta}}[0,\IF)
\Psi\LT(\frac{\sqrt{2}(\delta u+c)}{\sigma\sqrt{\delta}}\RT),
\end{align*}
where
\BQNY
\widetilde{\mathcal{P}}^{\frac{c^2}{\sigma^2\delta}}[0,\IF)=\E{\sup_{t\in [0, \IF)}
e^{\sqrt{2}B(t)-2t+\frac{2c}{\sigma\sqrt{\delta}}\sqrt{t}}}.
\EQNY
  The asymptotic result decreases when $\delta$ increases. Since $\widetilde{\mathcal{P}}^{\frac{c^2}{\sigma^2\delta}}[0,\IF)$ is a decreasing function of $\delta>0$  and $\Psi\LT(\frac{\sqrt{2}(\delta u+c)}{\sigma\sqrt{\delta}}\RT)$ is decreasing when $\delta$ increases and  $u\geq \frac{c}{\sqrt{2}\delta}$, the asymptotic of $\psi_{\gamma,\IF}^\delta(u)$ is also a decreasing function of $\delta$. The effect of $\delta$ is not an intuitionistic result from the original risk model.\\
Furthermore, comparing this result with it in \cite{Threshold2016,ParisianInfinite2018}, the scenario with tax is just $\frac{1}{1-\gamma}$ multiple of that without tax. Table \ref{tab2} is the simulated asymptotic results of $\psi_{\gamma,\IF}^\delta(u) $ in \netheo{TheomIF}.\\
ii) $\widehat{\mathcal{P}}^{\frac{c^2}{\sigma^2\delta}}[0,\IF)$ and  $\widetilde{\mathcal{P}}^{\frac{c^2}{\sigma^2\delta}}[0,\IF)$ can be considered as the generalised Piterbarg constants, see e.g., \cite{Pit96,MR1993262,debicki2011extremes,GeneralPit16}for various properties including the finiteness of related constants.\\
iii) We here interpret that the analysis of $\psi_{\gamma,\IF}^\delta(u) $ for the case $\delta<0$ is meaningless. We have
\BQNY
\sup_{t\in[0,\IF)}(u-U_\Ga^\delta(t))\geq \sup_{t\in[0,\IF)} \LT(\sigma\int_{0}^{t}e^{-\delta v}d B(v)-c\int_{0}^t e^{-\delta v}d v\RT)=\IF,\quad a.s.,
\EQNY
where $Var\LT(\sigma\int_{0}^{t}e^{-\delta v}d B(v)\RT)=\frac{\sigma^2}{2\delta}(1-e^{-2\delta t})$.
\end{remarks}
%\BK
%We have for $\Ga=0, \delta >0$ and $T=\IF$
%\BQNY
%\psi_{0,\IF}^\delta \sim
%\widehat{\mathcal{P}}^{1,\frac{c^2}{\sigma^2\delta}}[0,\IF)\Psi(M_u),
%\EQNY
%where $M_u$ and $\widehat{\mathcal{P}}^{1,\frac{c^2}{\sigma^2\delta}}[0,\IF)$ are the same as that of \netheo{TheomIF}.
%\EK
\begin{table}[H]
\centering
\caption{The simulated asymptotic results of $\psi_{\gamma,\IF}^\delta$}\label{tab2}
\begin{tabular}{p{1cm}p{1cm}p{1cm}p{1.5cm}p{1cm}p{4cm}}
\hline
$u$   & $c$   & $\sigma$ & $\delta$   & $\gamma$ & asymptotic results \\
\hline
5   & 0.1 & 1      & 0.05   &0.1    & 0.0467 \\
\hline
5   & 0.1 & 1      & 0.05   &0.2    & 0.0526 \\
\hline
5   & 0.1 & 1      & 0.07   &0.1    & 0.0256 \\
\hline
5   & 0.1 & 1      & 0.07   &0.2    & 0.0288 \\
\hline
5   & 0.1 & 1      & 0.1    &0.1   & 0.0113 \\
\hline
5   & 0.1 & 1      & 0.1    &0.2   & 0.0128 \\
\hline
4   & 0.1 & 1      & 0.1    &0.1   & 0.0378 \\
\hline
4   & 0.1 & 1      & 0.1    &0.2   & 0.0425 \\
\hline
\end{tabular}
%\\*we always set $u=5, c=0.1,\sigma=1, \lambda=600$.
\end{table}

We present below the approximation of the conditional passage time $\tau(u)|\tau(u)\leq T$ with $\tau(u)$ defined 	in \eqref{tu}.
\BT\label{ruintime}
For $T\in(0,\IF), \delta\in(-\IF,0)\cup(0,\IF)$ and $x>0$,  we have as $u\rw\IF$
\BQNY
\pk{u^2(T-\tau(u))>x|(\tau(u)\leq T)} \sim
\exp\LT(-\frac{\sigma^2e^{-2\delta T}x}{2a^2}\RT).
\EQNY
For $T=\IF, \delta>0$ and $x\in(-\frac{c}{\sigma\sqrt{\delta}},\IF)$, we have
\BQNY
\pk{u^{2}\LT(e^{-2\delta\tau(u)}-\LT(\frac{c}{\delta u+c}\RT)^2\RT)\leq x\big| \tau(u)<\IF}\sim \frac{\widehat{\mathcal{P}}^{\frac{c^2}{\sigma^2\delta}}[0,x+\frac{c}{\sigma\sqrt{\delta}}]}
{\widehat{\mathcal{P}}^{\frac{c^2}{\sigma^2\delta}}[0,\IF)},
\EQNY
as $u\rightarrow\IF$.
\ET
\begin{remark}
When $\gamma=0$, the result of the scenario $T=\IF$ corresponds to that in \cite{ParisianInfinite2018}.
\end{remark}
\section{Proofs}
Before giving the proofs of our main theorems, we need to introduce some notation which play a pivotal role in the proofs, starting with
\BQN\label{PS}
\mathcal{P}^a[0, S]=\E{ \sup_{t\in[0,S]}e^{\sqrt{2}B(t)-(1+a)\abs{t}}}\in (0,\IF),
\EQN
where $S, a$ are positive constants and
\BQNY
\mathcal{P}^a[0,\IF):=\underset{S\rw\IF}\lim\mathcal{P}^a[0,S]=1+\frac{1}{a}
\EQNY
where is known, see e.g., \cite{Pit96} or \cite{MR1993262}. \\
Moreover, Pickands constant defined by
$$\mathcal{H}[0,S]=\E{\sup_{t\in[0,S]} e^{\sqrt{2}B(t)-t}}\in(0,\IF).$$
It is known that $\lim_{T\rightarrow\infty} \frac{1}{T}\mathcal{H}[0,T]=1$, see  \cite{Pic1969, Berman82, Pit96,DebickiRol02, Mandjes2007,DMandjes2011,DiekerY} for various properties of $\mathcal{P}^a[0, S], \mathcal{H}[0,S]$ and its generalizations.\\
\prooftheo{TheomOT}
We define for any $\gamma \in (0,1)$ the random process $Z$ by
\BQN\label{ZZ1}
Z(s,t):=\frac{\sigma \int_0^t e^{-\delta v}dB(v)-c \int_0^t e^{-\delta v}dv-\Ga\left(\sigma \int_0^s e^{-\delta v}dB(v)-c \int_0^s e^{-\delta v}dv\right)}{1+\gamma(e^{-\delta s}-1)}, \quad s,t\ge 0,
\EQN
which is crucial for our analysis, then for any $u$ positive
\BQNY
\psi_{\gamma,T}^\delta(u)=\pk{\sup_{0\leq s\leq t\leq T}Z(s,t)>u}.
\EQNY
Define next the mean function of $Z(s,t)$
\BQNY
m(s,t):=\E{Z(s,t)}&=&\frac{1}{1+\gamma(e^{-\delta s}-1)}\left(-c\int_0^t e^{-\delta v}dv+\gamma c\int_0^s e^{-\delta v}dv\right)\\
&=&\frac{1}{1+\gamma(e^{-\delta s}-1)}\left(\frac{\gamma c}{\delta}\LT(1- e^{-\delta s}\RT)-\frac{c}{\delta}\LT(1- e^{-\delta t}\RT)\right)
\EQNY
and its variance function
\begin{align}
V^2_Z(s,t):&= \E{Z(s,t)-\E{Z(s,t)}}^2\label{VZ11}\\
&=\frac{1}{(1-\gamma+\gamma e^{-\delta s})^2}\E{\Bigl(\sigma\int_0^t e^{-\delta v}dB(v)-\gamma \sigma\int_0^s e^{-\delta v}dB(v)\Bigr)^2}\nonumber\\
&=\frac{\sigma^2}{2 \delta(1-\gamma+\gamma e^{-\delta s})^2}\LT(\LT(1-e^{-2\delta t}\RT)-\gamma(2-\gamma)\LT(1-e^{-2\delta s}\RT)\RT).\nonumber
\end{align}
By \nelem{VZ}, $V_Z^2(s,t)$ attains the unique maximum at the point$(0,T)$. We give the asymptotic expansion of the standard deviation function $V^2_Z(s,t)$ at this point. It follows that
\BQN\label{eq:VZ}
V_Z(s,t)=
a\LT(1-\frac{\Ga\sigma^2}{2a^2}(1-\gamma+e^{-2\delta T})s-\frac{\sigma^2e^{-2\delta T}}{2a^2}(T-t)\RT)+o((T-t)+s)
\EQN
as $(s,t)\rw(0,T)$; hence there exists a positive constant $\theta>0$ such that
\BQN\label{eq:det1}
  \abs{t-T-\gamma e^{\delta T}s}\le \mathbb{C} (\xx{V_Z}(0,T)-\xx{V_Z}(s,t))
\EQN
uniformly in $B_\theta:=\{(s,t): (s,t)\in[0,\theta]\times[T-\theta,T]\}$.
Besides, $m(0,T)=\frac{c}{\delta}(e^{-\delta T}-1)$.
Next, we study the asymptotic of the supremum of the Gaussian random  field defined on $B_\theta$.
Set below
\BQNY
\nu_u(s,t)=\frac{u-m(s,t)}{V_Z(s,t)},\ W(s,t)=Z(s,t)-m(s,t).
\EQNY
Setting $V^2_W(s,t)=Var(W(s,t))$ and $\overline{W}(s,t)=\frac{W(s,t)}{V_W(s,t)}$, it is clear that $V_W=V_Z$, for any $u>0$
\begin{align}\label{piZ}
\Pi(u)\le\pk{\sup_{(s,t)\in B}Z(s,t)>u}
\le \Pi(u)+\Pi_1(u),
\end{align}
where
\BQNY
&&\Pi(u)=\pk{\sup_{(s,t)\in B_\theta}\overline{W}(s,t)\frac{\nu_u(0,T)}{\nu_u(s,t)}>\nu_u(0,T)},\\
&&\Pi_1(u)=\pk{\sup_{(s,t)\in B\setminus B_\theta}\overline{W}(s,t)\frac{\nu_u(0,T)}{\nu_u(s,t)}>\nu_u(0,T)}.
\EQNY
Since
\begin{align}\label{bibi}
\frac{\nu_u(0,T)}{\nu_u(s,t)}&=1-\frac{V_W(0,T)-V_W(s,t)}{V_W(0,T)}+\frac{[m(s,t)-m(0,T)]V_W(s,t)}{(u-m(s,t))V_W(0,T)}\nonumber\\
&=1-\frac{V_W(0,T)-V_W(s,t)}{V_W(0,T)}+\frac{ce^{-\delta T}[(T-t)+\gamma e^{\delta T} s]V_W(s,t)+o((T-t)+s)}{(u-m(s,t))V_W(0,T)},
\end{align}
as $(s,t)\rightarrow(0,T)$, we have, in view of \eqref{eq:det1}, for any $\vp\in(0,1)$,  and sufficiently large $u$
\BQNY
1-\frac{V_W(0,T)-V_W(s,t)}{V_W(0,T)}\le \frac{\nu_u(0,T)}{\nu_u(s,t)}\le 1-(1-\vp)\frac{V_W(0,T)-V_W(s,t)}{V_W(0,T)}
\EQNY
uniformly in $(s,t)\in B_\theta$.
Consequently
\BQN
\pk{\sup_{(s,t)\in B_\theta}W_0(s,t)>\nu_u(0,T)}\le \Pi(u)\le \pk{\sup_{(s,t)\in B_\theta}W_\vp(s,t)>\nu_u(0,T)},
\EQN
where the random field $\{W_\vp(s,t), s,t\ge0\}$ is defined as
\BQNY
W_\vp(s,t):=\overline{W}(s,t)\left(1-(1-\vp)\frac{V_W(0,T)-V_W(s,t)}{V_W(0,T)}\right),\ \  \vp\in[0,1).
\EQNY
Direct calculations show that the standard deviation function $\sigma_{W_{\vp}}(s,t):=\sqrt{\E{(W_\vp(s,t))^2}}$ attains its unique maximum over $B_\theta$ at $(0,T)$ with $\sigma_{W_{\vp}}(0,T)=1$.
Thus, in the light of \eqref{eq:VZ}, we have
\BQN\label{eq:SWV}
\sigma_{W_\varepsilon}(s,t)=1-(1-\vp)\left(\frac{\Ga\sigma^2(1-\Ga+e^{-2\delta T})}{2a^2}s+\frac{\sigma^2e^{-2\delta T}}{2a^2}(T-t)\right)(1+o(1)),
\EQN
as $(s,t)\rw(0,T)$. Furthermore, it follows that
\BQN\label{cove1}
1-Cov(W_\vp(s,t),W_\vp(s',t'))=\frac{\sigma^2}{2a^2}\left(e^{-2\delta T}\mid t-t'\mid +\Ga^{2}\mid s-s'\mid\right)(1+o(1))
\EQN
as $(s,t), (s',t')\rw(0,T)$. In addition, we obtain
\BQNY
\E{(W_\vp(s,t)-W_\vp(s',t'))^2}\le \mathbb{C}(2e^{-2\delta T}|t-t'|+2\Ga^2|s-s'|)
\EQNY
for $(s,t), (s',t')\in B_\theta$. Consequently, by Theorem 8.2 of \cite{Pit96}
\BQN
\pk{\sup_{(s,t)\in B_\theta}W_\vp(s,t)>\nu_u(0,T)}\sim
\frac{2\LT(1+e^{-2\delta T}\RT)}{1-\gamma+e^{-2\delta T}}\Psi\left(\frac{u+\frac{c}{\delta}(1-e^{-\delta T})}{a}\right)\label{result}
\EQN
 as $u\rw\infty,\ \vp\rw 0$. Thus we obtain the asymptotic upper bound for $\Pi(u)$ on the set $B_{\theta}$. The asymptotic lower bound can be derived using the same arguments. In order to complete the proof we need to show further that
\BQN
\Pi_1(u)=o(\Pi(u))\ \ \text{as}\ u\rw\IF.
\label{WB}
\EQN
In the light of \eqref{bibi} for all $u$ sufficiently large
\BQNY
\sup_{(s,t)\in B\setminus B_\theta} Var\left(\overline{W}(s,t)\frac{\nu_u(0,T)}{\nu_u(s,t)}\right)\le (\rho(\theta))^2<1,
\EQNY
where $\rho(\theta)$ is a positive function in $\theta$ which exists due to the continuity of $V_W(s,t)$ in $B$.
Therefore, a direct application of Borell-TIS inequality as in \cite{AdlerTaylor}  implies
\BQNY
\pk{\sup_{(s,t)\in B \setminus B_\theta}\overline{W}(s,t)\frac{\nu_u(0,T)}{\nu_u(s,t)}>\nu_u(0,T)}\le 2\Psi\LT(\frac{\nu_u(0,T)-b}{\rho(\theta)}\RT)=o(\Pi(u)), \quad u\rw\IF,
\EQNY
where $b=\sup_{(s,t)\in B \setminus B_\theta}\E{\overline{W}(s,t)\frac{\nu_u(0,T)}{\nu_u(s,t)}}<\IF.$\\
Consequently, Eq. (\ref{WB}) is established.
\QED
\noindent
\prooftheo{TheomIF}
We have
\BQNY
\psi_{\gamma,\IF}^\delta(u)=\pk{\sup_{0<t\leq s\leq 1}\frac{Z(s,t)}{G_u(s,t)}>1},
\EQNY
 where
\begin{align}
Z(s,t)&=\sigma\int_0^{-\frac{1}{2\delta}lnt} e^{-\delta v}dB(v)-\Ga\sigma\int_0^{-\frac{1}{2\delta}lns} e^{-\delta v}dB(v),\label{ZZ2}\\
G_u(s,t)&=u+\Ga(s^{\frac{1}{2}}-1)u+c\int_0^{-\frac{1}{2\delta}lnt} e^{-\delta v}dv-c\Ga\int_0^{-\frac{1}{2\delta}lns} e^{-\delta v}dv\nonumber\\
&=u-\Ga \LT(u+\frac{c}{\delta}\RT)(1-s^{\frac{1}{2}})+\frac{c}{\delta}\LT(1-t^{\frac{1}{2}}\RT).\label{GU2}
\end{align}
The variance function of $Z(s,t)$ is given by
\begin{align}
V_Z^2(s,t)&=Var\left(\sigma\int_0^{-\frac{1}{2\delta}lnt} e^{-\delta v}dB(v)-\Ga\sigma\int_0^{-\frac{1}{2\delta}lns} e^{-\delta v}dB(v)\right)\label{VZ22}\\
&=\frac{\sigma^2}{2\delta}((1-t)-\Ga(2-\Ga)(1-s)).\nonumber
\end{align}
Let $M_u(s,t)=\frac{G_u(s,t)}{V_Z(s,t)}$, then for $t_u$ in \nelem{VIF}
\BQNY
M_u:=M_u(1,t_u)=\frac{\sqrt{2}}{\sigma}\sqrt{\delta u^2+2cu}=\frac{\sqrt{2\delta}}{\sigma}u(1+o(1)),
\EQNY
as $u\rw\IF$.
\BQNY
\frac{M_u(s,t)}{M_u}-1=\frac{[G_u(s,t)V_Z(1,t_u)]^2-[G_u(1,t_u)V_Z(s,t)]^2}{V_Z(s,t)G_u(1,t_u)[G_u(s,t)V_Z(1,t_u)+V_Z(s,t)G_u(1,t_u)]}
\EQNY
Since
\begin{align*}
&[G_u(s,t)V_Z(1,t_u)]^2-[G_u(1,t_u)V_Z(s,t)]^2\\
=&\{(u+\frac{c}{\delta})[1-\Ga(1-\sqrt{s})]-\frac{c}{\delta}\sqrt{t}\}^2\frac{\sigma^2}{2\delta}(1-t_u)\\
&-\{(u+\frac{c}{\delta})-\frac{c}{\delta}\sqrt{t_u}\}^2\frac{\sigma^2}{2\delta}[1-t-\Ga(2-\Ga)(1-s)]\\
\sim&\frac{\sigma^2}{2\delta}u^2[(\sqrt{t}-\sqrt{t_u})^2+\Ga(1-\Ga)(1-s)],
\end{align*}
then
\BQN\label{MM1}
\frac{M_u(s,t)}{M_u}-1\sim \frac{1}{2}[(\sqrt{t}-\sqrt{t_u})^2+\Ga(1-\Ga)(1-s)].
\EQN
Now we rewrite
\BQNY
\psi_{\gamma,\IF}^\delta(u)=\pk{\sup_{0<t\leq s\leq 1}\frac{Z(s,t)}{G_u(s,t)}M_u>M_u}.
\EQNY
The correlation function of $Z(s,t)$ is
\BQNY
r(s,s',t,t')
%&=&\frac{\mathbb{E}\left[(\sigma\int_0^{-\frac{1}{2\delta}lnt} e^{-\delta v}dB(v)-\Ga\sigma\int_0^{-\frac{1}{2\delta}lns} e^{-\delta v}dB(v))(\sigma\int_0^{-\frac{1}{2\delta}lnt'} e^{-\delta v}dB(v)-\Ga\sigma\int_0^{-\frac{1}{2\delta}lns'} e^{-\delta v}dB(v))\right]} {\sqrt{\frac{\sigma^2}{2\delta}[1-t-\Ga(2-\Ga)(1-s)]}\sqrt{\frac{\sigma^2}{2\delta}[1-t'-\Ga(2-\Ga)(1-s')]}}\\
\sim 1-\frac{1}{2}|t-t'|-\frac{1}{2}\Ga^2|s-s'|
\EQNY
for $s,s'\rw1$, $t-t_u, t'-t_u\rw0$.\\
In addition, we obtain for some $\theta_0,\mathbb{C},$
\BQNY
\mathbb{E}\left(Z(s,t)-Z(s',t')\right)^2\leq\mathbb{C}(|t-t'|+\Ga^2|s-s'|)
\EQNY
for $ s,t\in[1-\theta_0,1]\times[0,\theta_0]$.\\
\def\TPi{\widetilde{\Pi}}
Note that for any small $\theta_1,\theta_2\in(0,1)$, set $B=\{(s,t):0<t\leq s\leq 1\}$ and $\Delta_\theta=[1-\theta_1,1]\times[0,\theta_2]$,
\BQNY
\Pi(u):=\pk{\sup_{(s,t)\in\Delta_\theta}\frac{Z(s,t)}{G_u(s,t)}M_u>M_u}\leq \psi_{\gamma,\IF}^\delta(u)\leq\Pi(u) +\TPi(u)
\EQNY
with
\BQNY
\TPi(u):=\pk{\sup_{t\in B\setminus\Delta_\theta}\frac{Z(s,t)}{G_u(s,t)}M_u>M_u}
\EQNY
 In the following, we shall focus on the asymptotics of $\Pi(u)$ as $u\to\IF$, and finally we show that
 \BQN\label{PIPI}
 \TPi(u)=o(\Pi(u)),\ \ \ u\to\IF.
 \EQN
We set $\delta_1(u)=\LT(\frac{(lnu)^q}{u}\RT)^2$, $\delta_2(u)=2\sqrt{t_u}\frac{(lnu)^q}{u}+\LT(\frac{(lnu)^q}{u}\RT)^2$ for some $q>1$ and
$$D_u=[1-\delta_1(u),1]\times[0,t_u+\delta_2(u)]=\LT[0,\LT(\frac{(lnu)^q}{u}\RT)^2\RT]\times\LT[0,\LT(\sqrt{t_u}+\frac{(lnu)^q}{u}\RT)^2\RT],\ \ \Theta_u=\Delta_\theta\setminus D_u.$$
Clearly, for $u$ large enough,
\BQNY
\Pi_1(u):=\pk{\sup_{t\in D_u}\frac{Z(s,t)}{G_u(s,t)}M_u>M_u}\leq \Pi(u)\leq \Pi_1(u)+\Pi_2(u),
\EQNY
with
\BQNY
\Pi_2(u)=\pk{\sup_{t\in \Theta_u}\frac{Z(s,t)}{G_u(s,t)}M_u>M_u}
\EQNY
Next, we give a tight upper bound  for $\Pi_2(u)$ which will finally imply that, for some small $\theta_1, \theta_2\in(0,1)$
\BQNY
\Pi_2(u)=o(\Pi_1(u)),\ \ \ u\to\IF.
\EQNY
By (\ref{MM1}) we have that for any small $\vn_o>0$, there exists some small $\theta_1, \theta_2>0$ such that
\BQNY
\frac{M_u}{M_u(s,t)}\leq\frac{1}{1+\frac{1-\vn_0}{2}[(\sqrt{t}-\sqrt{t_u})^2+\Ga(1-\Ga)(1-s)]}
\EQNY
holds for all $(s,t)\in \Delta_\theta$. Furthermore, for any $t\in\Theta_u$
\BQNY
1+\frac{1-\vn_0}{2}[(\sqrt{t}-\sqrt{t_u})^2+\Ga(1-\Ga)(1-s)]\geq 1+\frac{1-\vn_0}{2}\min(1, \Ga(1-\Ga))\LT(\frac{(lnu)^q}{u}\RT)^2
\EQNY
implying $\left(set \overline{Z}(s,t)=\frac{Z(s,t)}{V_Z(s,t)}\right)$
\BQN
\Pi_2(u)
&=&\pk{\sup_{(s,t)\in\Theta_u} \overline{Z}(s,t)\frac{M_u}{M_u(s,t)} >M_u}\nonumber\\
&\leq& \mathbb{Q}_1M_u^4\Psi\left(M_u \LT(1+\frac{1-\vn_0}{2}\min(1, \Ga(1-\Ga))\LT(\frac{(lnu)^q}{u}\RT)^2\RT)\right)
\EQN
holds for all $u$ large, with some constant $\mathbb{Q}_1>0$.
For any small $\vn>0$, define below
$$ \mathcal{B}_u^{\vn\pm}(\Delta):= \Bigl\{ \sup_{(s,t)\in \Delta}  \frac{\xi^{\pm}(s,t)}{[1+(\frac{1}{2}\pm\vn)(\sqrt{t}-\sqrt{t_u})^2][1+(\frac{1}{2}\Ga(1-\Ga)\pm\vn)(1-s)]}>M_u\Bigr\},\ \Delta\subset \R^2,$$
where $\{\xi^{\pm}(s,t), s, t\geq 0\}$ is a zero-mean stationary Gaussian field with continuous sample paths and correlation function
\BQNY
 \Ga^{\pm}_\xi(s,t)=\exp\LT(-\LT(\frac{{\Ga}^2}{2}\mp\vn\RT)|s|-\LT(\frac{1}{2}\mp\vn\RT)|t|\RT).
\EQNY
Next we analyse $\Pi_1(u)$, as $u\rw\IF$. By (\ref{MM1}), when $u$ large enough, for any $t\in D_u$,
\BQNY
&&\frac{M_u}{M_u(s,t)}\geq\frac{1}{[1+(\frac{1}{2}+\vn)(\sqrt{t}-\sqrt{t_u})^2]
[1+(\frac{1}{2}\Ga(1-\Ga)+\vn)(1-s)]},\\
&& \frac{M_u}{M_u(s,t)}\leq
\frac{1}{[1+(\frac{1}{2}-\vn)(\sqrt{t}-\sqrt{t_u})^2][1+(\frac{1}{2}\Ga(1-\Ga)-\vn)(1-s)]}.
\EQNY
Then
\BQNY
\mathbb{P}\{\mathcal{B}_u^{\vn+}(D_u)\}\leq\Pi_1(u)\leq\mathbb{P}\{\mathcal{B}_u^{\vn-}(D_u)\},\ u\rw\IF.
\EQNY
Thus we just need establish the asymptotic behavior of $\pi^{\pm}(u):=\mathbb{P}\{\mathcal{B}_u^{\vn\pm}(D_u)\}$, then according to the continuous of the results which can be seen from the following calculation, setting $\vn\rw 0$, we will gain the precision estimates of $\Pi_1(u)$. Below we mainly show the calculation of $\pi^{+}(u)$.

$$D_u=[1-\delta_1(u),1]\times[0,t_u+\delta_2(u)]=\LT[0,\LT(\frac{(lnu)^q}{u}\RT)^2\RT]\times\LT[0,\LT(\sqrt{t_u}+\frac{(lnu)^q}{u}\RT)^2\RT]$$
For any positive constant $S_1,S_2$, define
\BQNY
&&\Delta^1_k=1-u^{-2}[(k+1)S_1,kS_1],\ \ k=0,1,2,3\ldots\\
&&\Delta^2_{-1}=[0,t_u],\ \ \Delta^2_k=[(\sqrt{t_u}+\sqrt{kS_2}u^{-1})^2,(\sqrt{t_u}+\sqrt{(k+1)S_2}u^{-1})^2],\ \ k=0,1,,3\ldots
\EQNY
and let further for $u>0$
\BQNY
h_1(u)=\lfloor S_1^{-1}(\ln u)^{2q}\rfloor+1,\ \
h_2(u)=\LT\lfloor \frac{\LT(\sqrt{t_u}+\frac{(\ln u)^q}{u}\RT)^2-t_u}{(\sqrt{t_u}+\sqrt{(k+1)S_2}u^{-1})^2
-(\sqrt{t_u}+\sqrt{kS_2}u^{-1})^2}\RT\rfloor+1,
\EQNY
where $h_1(u),h_2(u)\rw\IF$, as $u\rw\IF$.
By Bonferroni's inequality we have
\BQN\label{piz1}
\pi_1^+(u):=\pk{\mathcal{B}^{\vn+}_u(\Delta^1_0\times(\Delta^2_{-1}\cup\Delta^2_0))}
\leq\pi^{+}(u)\leq\pi_1^+(u)+\pi_2^+(u),
\EQN
where
\begin{align*}
\pi_2^+(u)=&\sum_{k_1=1}^{h_1(u)}\pk{\mathcal{B}^{\vn+}_u(\Delta^1_{k_1}\times(\Delta^2_{-1}\cup\Delta^2_0))}
+\sum_{k_2=1}^{h_2(u)}\pk{\mathcal{B}^{\vn+}_u(\Delta^1_0\times\Delta^2_{k_2})}\nonumber\\
&+\sum_{k_1=1}^{h_1(u)}\sum_{k_2=1}^{h_2(u)}\pk{\mathcal{B}^{\vn+}_u(\Delta^1_{k_1}\times\Delta^2_{k_2})}\\
=:& I_1+I_2+I_3.\nonumber
\end{align*}
By \textbf{\nelem{PST}}, we get as $u\rw\IF$
\BQN
\pi_1^+(u)\sim\mathcal{P}^{\frac{\Ga(1-\Ga)+2\vn}{\Ga^2-2\vn}}[0,(\Ga^2-2\vn)\frac{\delta}{\sigma^2}S_1]
\widehat{\mathcal{P}}^{\frac{1+2\vn}{1-2\vn},\frac{c^2}{\sigma^2\delta}(1-2\vn)}[0,(1-2\vn)\frac{\delta}{\sigma^2}(\frac{c}{\delta}+\sqrt{S_2})^2]\Psi(M_u).
\EQN
Next we calculate the required asymptotic bounds for $I_1(u)$, $I_2(u)$ and $I_3(u)$ and show that as $u\rw\IF,\ S_i\rw\IF,\ i=1,2 $
\BQN\label{II}
I_1(u)=o(\pi_1^+(u)),I_2(u)=o(\pi_1^+(u)),I_3(u)=o(\pi_1^+(u)).
\EQN
By \textbf{\nelem{PST}}, we drive that
\begin{align*}
I_1(u)=&\sum_{k_1=1}^{h_1(u)}\pk{\mathcal{B}^{\vn+}_u(\Delta^1_{k_1}\times(\Delta^2_{-1}\cup\Delta^2_0))}\\
\leq&\sum_{k_1=1}^{h_1(u)}\pk{\sup_{(s,t)\in(\Delta^1_{k_1}\times(\Delta^2_{-1}\cup\Delta^2_0))}
\frac{\xi^{+}(s,t)}{1+(\frac{1}{2}+\vn)(\sqrt{t}-\sqrt{t_u})^2}>M_u[1+(\frac{1}{2}\Ga(1-\Ga)+\vn)u^{-2}k_1S_1]}\\
\leq&\mathcal{H}[0,(\Ga^2-2\vn)\frac{\delta}{\sigma^2}S_1]
\widehat{\mathcal{P}}^{\frac{1+2\vn}{1-2\vn},\frac{c^2}{\sigma^2\delta}(1-2\vn)}[0,(1-2\vn)\frac{\delta}{\sigma^2}(\frac{c}{\delta}+\sqrt{S_2})^2]\frac{1}{\sqrt{2\pi}}\\
&\times\sum_{k_1=1}^{h_1(u)}\frac{1}{M_u[1+(\frac{1}{2}\Ga(1-\Ga)+\vn)u^{-2}k_1S_1]}
\exp\LT(-\frac{M_u^2[1+(\frac{1}{2}\Ga(1-\Ga)+\vn)u^{-2}k_1S_1]^2}{2}\RT)(1+o(1))\\
=&\mathcal{H}[0,(\Ga^2-2\vn)\frac{\delta}{\sigma^2}S_1]
\widehat{\mathcal{P}}^{\frac{1+2\vn}{1-2\vn},\frac{c^2}{\sigma^2\delta}(1-2\vn)}[0,(1-2\vn)\frac{\delta}{\sigma^2}(\frac{c}{\delta}+\sqrt{S_2})^2]\\
&\times\Psi(M_u)\sum_{k_1=1}^{h_1(u)}\exp\LT(-\frac{\delta}{\sigma^2}(\Ga(1-\Ga)+2\vn)k_1S_1\RT)(1+o(1))
\end{align*}
as $u\rw\IF$.
Similarly
\begin{align*}
I_2(u)\leq& \mathcal{P}^{\frac{\Ga(1-\Ga)+2\vn}{\Ga^2-2\vn}}[0,(\Ga^2-2\vn)\frac{\delta}{\sigma^2}S_1]
\widehat{\mathcal{P}}^{0,\frac{c^2}{\sigma^2\delta}(1-2\vn)}
[0,(1-2\vn)\frac{\delta}{\sigma^2}(\frac{c}{\delta}+\sqrt{S_2})^2]\\
&\times\Psi(M_u)\sum_{k_2=1}^{h_2(u)}\exp\LT(-\frac{\delta}{\sigma^2}(1+2\vn)k_2S_2\RT)(1+o(1)),
\end{align*}
as $u\rw\IF$. Moreover,
\begin{align*}
I_3(u)=&\sum_{k_1=1}^{h_1(u)}\sum_{k_2=1}^{h_2(u)}\pk{\mathcal{B}^{\vn+}_u(\Delta^1_{k_1}\times\Delta^2_{k_2})}\\
\leq&\sum_{k_1=1}^{h_1(u)}\sum_{k_2=1}^{h_2(u)}\pk{\sup_{(s,t)\in(\Delta^1_{k_1}\times(\Delta^2_{-1}\cup\Delta^2_0))}\xi^{+}(s,t)
>M_u[1+(\frac{1}{2}\Ga(1-\Ga)+\vn)u^{-2}k_1S_1+(\frac{1}{2}+\vn)u^{-2}k_2S_2]}\\
=&\mathcal{H}[0,(\Ga^2-2\vn)\frac{\delta}{\sigma^2}S_1]
\widehat{\mathcal{P}}^{0,\frac{c^2}{\sigma^2\delta}(1-2\vn)}[0,(1-2\vn)\frac{\delta}{\sigma^2}(\frac{c}{\delta}+\sqrt{S_2})^2]\\
&\times\Psi(M_u)\sum_{k_1=1}^{h_1(u)}\sum_{k_2=1}^{h_2(u)}\exp\LT(-\frac{\delta}{\sigma^2}[(\Ga(1-\Ga)+2\vn)k_1S_1+(1+2\vn)k_2S_2]\RT)
(1+o(1)),
\end{align*}
as $u\rw\IF$. By letting $S_1,S_2\rw\IF$, (\ref{II}) is proved.\\
Thus we finish the calculation of $\pi^+(u)$.\\
For $\pi^-(u)$, we just need notice that the equations (\ref{piz1}) are replaced by
\BQNY
\pi_1^-(u):=\pk{\mathcal{B}^{\vn-}_u(\Delta^1_0\times(\Delta^2_{-1}\cup\Delta^2_0))}\leq\pi^{-}(u)\leq\pi_1^-(u)+\pi_2^-(u)
\EQNY
and
\begin{align*}
\pi_2^-(u)=&\sum_{k_1=1}^{h_1(u)}\pk{\mathcal{B}^{\vn-}_u(\Delta^1_{k_1}\times(\Delta^2_{-1}\cup\Delta^2_0))}
+\sum_{k_2=1}^{h_2(u)}\pk{\mathcal{B}^{\vn-}_u(\Delta^1_0\times\Delta^2_{k_2})}\\
&+\sum_{k_1=1}^{h_1(u)}\sum_{k_2=1}^{h_2(u)}\pk{\mathcal{B}^{\vn-}_u(\Delta^1_{k_1}\times\Delta^2_{k_2})}\\
=&: J_1+J_2+J_3.
\end{align*}
Then by \textbf{\nelem{PST}}, we get
\BQNY
\pi_1^-(u)\sim\mathcal{P}^{\frac{\Ga(1-\Ga)-2\vn}{\Ga^2+2\vn}}[0,(\Ga^2+2\vn)\frac{\delta}{\sigma^2}S_1]
\widehat{\mathcal{P}}^{\frac{1-2\vn}{1+2\vn},\frac{c^2}{\sigma^2\delta}(1+2\vn)}
[0,(1+2\vn)\frac{\delta}{\sigma^2}(\frac{c}{\delta}+\sqrt{S_2})^2]\Psi(M_u),
\EQNY
and similarly
\BQNY
J_1(u)=o(\pi_1^-(u)),J_2(u)=o(\pi_1^-(u)),J_3(u)=o(\pi_1^-(u)),\ as\  u\rw\IF,\ S_i\rw\IF,\ i=1,2.
\EQNY
Thus,
$$\lim_{\vn\rw 0} \pi^+(u)=\lim_{\vn\rw 0} \pi^-(u)=\Pi_1(u).$$
Finally, by (\ref{MM1}) and \textbf{\nelem{VIF}}, we can choose some small $\theta_1, \theta_2>0$ so that for any $u$ sufficiently large
\BQN
\sup_{(s,t)\in B\setminus\Delta_\theta}Var\LT\{\frac{Z(s,t)}{G_u(s,t)}M_u\RT\}\leq\sup_{(s,t)\in B\setminus\Delta_\theta}\LT(\frac{M_u}{M_u(s,t)}\RT)^2\leq (\rho(\theta_1, \theta_2))^2< 1
\EQN
where $\rho(\theta_1, \theta_2)$ is a positive function in $\theta_1$ and $\theta_2$ which exists due to the continuity of $M_u(s,t)$ in $B$. Additionally, by the almost surely continuity of random field, we have, for some constant $b>0$

\BQNY
\pk{\sup_{t\in B\setminus\Delta_\theta}\frac{Z(s,t)}{G_u(s,t)}M_u>b}\leq \frac{1}{2}.
\EQNY
Therefore, a direct application of the Borell inequality (e.g., \textbf{Theorem D.1} of \cite{Pit96})
 implies
 \BQNY
 \TPi(u)=\pk{\sup_{t\in B\setminus\Delta_\theta}\frac{Z(s,t)}{G_u(s,t)}M_u>M_u}\leq 2\Psi\LT(\frac{M_u-b}{\rho(\theta_1, \theta_2)}\RT)=o(\Pi(u))\
 \ as\ \ u\rw \IF.
 \EQNY
 Consequently, Eq. (\ref{PIPI}) is established, and thus the proof is complete.
 \QED
 \prooftheo{ruintime} For $T\in(0,\IF)$, first note that
\BQN\label{eq:PP}
\pk{u^2(T-\tau(u))>x|(\tau(u)\leq T)}=\frac{\pk{\sup_{0\leq s\leq t\leq T_x(u)}Z(s,t)>u}}{\pk{\sup_{0\leq s\leq t\leq T}Z(s,t)>u}}
\EQN
for any $x,u>0,$ where $T_x(u)=T-xu^{-2}$ and $Z(s,t)$ is the same as in \eqref{ZZ1}. Next for any $x,u>0$, we follow the similar argumentation as in the proof \netheo{TheomOT}
$$\pk{\sup_{0\leq s\leq t\leq T_x(u)}Z(s,t)>u}\sim \mathcal{P}^{\frac{1-\gamma+e^{-2\delta T}}{\gamma}}\lim_{S\rw\IF}\mathcal{P}^{1}\LT[\frac{\sigma^2e^{-2\delta T}}{2a^2}x,\frac{\sigma^2e^{-2\delta T}}{2a^2}S\RT]\Psi\LT(\frac{\sqrt{2}}{\sigma}\sqrt{\delta u^2+2cu}\RT),$$
as $u\rw\IF$.
By Remark 2.3 of \cite{Threshold2016},  we have
\begin{align*}
\lim_{S\rw\IF}\mathcal{P}^{1}\LT[\frac{\sigma^2e^{-2\delta T}}{2a^2}x,\frac{\sigma^2e^{-2\delta T}}{2a^2}S\RT]
&=\lim_{S\rw\IF}\E{ \sup_{t\in\LT[\frac{\sigma^2e^{-2\delta T}}{2a^2}x,\frac{\sigma^2e^{-2\delta T}}{2a^2}S\RT]}e^{\sqrt{2}B(t)-2t}}\\
&=\lim_{S\rw\IF}\E{ \sup_{t\in\LT[0,\frac{\sigma^2e^{-2\delta T}}{2a^2}(S-x)\RT]}e^{\sqrt{2}B(t)-2t-\frac{\sigma^2e^{-2\delta T}}{2a^2}x}}\\
&=e^{-\frac{\sigma^2e^{-2\delta T}}{2a^2}x}\lim_{S\rw\IF}\E{ \sup_{t\in\LT[0,\frac{\sigma^2e^{-2\delta T}}{2a^2}(S-x)\RT]}e^{\sqrt{2}B(t)-2t}}\\
&=e^{-\frac{\sigma^2e^{-2\delta T}}{2a^2}x}\mathcal{P}^{1}[0,\IF)\\
&=2e^{-\frac{\sigma^2e^{-2\delta T}}{2a^2}x}.
\end{align*}
Therefore, we conclude by an application of Theorem \ref{TheomOT} and equation(\ref{eq:PP}) that
$$\pk{u^2(T-\tau_1(u))>x|(\tau_1(u)\leq T)}\rightarrow \exp\LT(-\frac{\sigma^2e^{-2\delta T}x}{2a^2}\RT), $$
as $u\rightarrow \IF$, for any $x>0$.\\
For the case $T=\IF$, we have
\BQNY
\pk{u^{2}\LT(e^{-2\delta\tau_u}-\LT(\frac{c}{\delta u+c}\RT)^2\RT)\leq x\big| \tau_u<\IF}
=\frac{\pk{\underset{0<t\leq t_u+u^{-2}x}{\sup_{0<t\leq s\leq 1}}\frac{Z(s,t)}{G_u(s,t)}>1}}{\pk{\sup_{0<t\leq s\leq 1}\frac{Z(s,t)}{G_u(s,t)}>1}}
\EQNY
where $Z(s,t)$ and $G_u(s,t)$ are the same as in \eqref{ZZ2} and \eqref{GU2}.
We follow the similar argumentation as in the proof \netheo{TheomIF}
\BQNY
\pk{\underset{0<t\leq t_u+u^{-2}x}{\sup_{0<t\leq s\leq 1}}\frac{Z(s,t)}{G_u(s,t)}>1}\sim \frac{1}{1-\gamma}\widehat{\mathcal{P}}^{1,\frac{c^2}{\sigma^2\delta}}[0,x+\frac{c}{\sigma\sqrt{\delta}}]
\Psi\LT(\frac{\sqrt{2}}{\sigma}\sqrt{\delta u^2+2cu}\RT),
\EQNY
as $u\rw\IF$.
Thus we get the results by an application of Theorem \ref{TheomIF}.
\QED
\section{Appendix}
\noindent
Here  we give several Lemmas which are used in the proofs.

\BEL\label{VZ}
The variance function $V^2_Z(s,t)$ in \eqref{VZ11} attains its unique global maximum over set $B:=\{(s,t):0\leq s\leq t<T\}$ at $(s_0,t_0)$, with $s_0=0$ and
$t_0= T$. Further
$
a^2:=V_Z^2(0,T)=\frac{\sigma^2}{2\delta}(1-e^{-2\delta T}).
$
\EEL
\noindent
\prooflem{VZ}
It is obvious that $t_0=T$, then
\begin{align*}
\frac{\partial V_Z^2(s,T)}{\partial s}&=\frac{\gamma \sigma^2 e^{-\delta s}}{(1-\gamma+\gamma e^{-\delta s})^3}\left((1-e^{-2\delta T})-\gamma(2-\gamma)\LT(1-e^{-2\delta s}\RT)\right)-\frac{\gamma (2-\gamma) \sigma^2 e^{-2\delta s}}{(1-\gamma+\gamma e^{-\delta s})^2}\\
%&=&\frac{\gamma \sigma^2 e^{-\delta s}}{(1-\gamma+\gamma e^{-\delta s})^3}\left[
%(1-e^{-2\delta T})-\gamma(2-\gamma)+\gamma(2-\gamma)e^{-2\delta s}-(2-\gamma)e^{-\delta s}(1-\gamma+\gamma e^{-\delta s})\right]\\
&=\frac{\gamma \sigma^2 e^{-\delta s}}{(1-\gamma+\gamma e^{-\delta s})^3}\left(
(1-e^{-2\delta T})-\gamma(2-\gamma)-(1-\gamma)(2-\gamma) e^{-\delta s}
\right),
\end{align*}
we have when $\delta<0$, $V_Z^2(s,T)$ attains the maximum only at $s=0$ and when $\delta>0$, $V_Z^2(s,T)$ attains the maximum only at $s=0$ and $s=T$, since
\BQNY
V_Z^2(0,T)=\frac{\sigma^2}{2\delta}(1-e^{-2\delta T}),\ \  \
V_Z^2(T,T)=\frac{\sigma^2(1-\gamma)^2}{2\delta(1-\gamma+\gamma e^{-\delta T})}(1-e^{-2\delta T}),
\EQNY
and $
\frac{(1-\gamma)^2}{(1-\gamma+\gamma e^{-\delta T})}<1,
$
hence the claim follows. \QED

\BEL\label{VIF}
Let $F_u(s,t)=\frac{V_Z(s,t)}{G_u(s,t)}$ with  $G_u(s,t)$ in \eqref{GU2} and $V_Z(s,t)$ in \eqref{VZ22}. Then for $u$ sufficiently large, the function $F_u(s,t), 0<t \leq s\leq 1 $ attains its maximum at the unique point $(s,t)=(1,t_u)$, where
$$t_u=\LT(\frac{c}{\delta u+c}\RT)^2.$$
\EEL
\noindent\prooflem{VIF} Note that
\BQNY
\frac{\partial V_Z^2(s,t)}{\partial t}=2V_Z(s,t)\frac{\partial V_Z(s,t)}{ \partial t}=-\frac{\sigma^2}{2\delta}, \ \ \frac{\partial V_Z^2(s,t)}{\partial s}=2V_Z(s,t)\frac{\partial V_Z(s,t)}{\partial s}=\frac{\sigma^2}{2\delta}\Ga(2-\Ga).
\EQNY
Then we have
\BQNY
\frac{\partial F_u(s,t)}{\partial t}&=&\frac{\partial V_Z(s,t)}{\partial t}\cdot\frac{1}{G_u(s,t)}-\frac{V_Z(s,t)}{G_u^2(s,t)}\LT(-\frac{c}{2\delta}t^{-\frac{1}{2}}\RT)\\
&=&\frac{1}{2G_u^2(s,t) V_z(s,t)}\left(\frac{\partial V_Z^2(s,t)}{\partial t}G_u(s,t)+V_Z^2(s,t)\frac{ct^{-\frac{1}{2}}}{\delta}\right)\\
&=&\frac{\sigma^2 t^{-1/2}}{4\delta G_u^2(s,t) V_Z(s,t)}\left\{[(1-\Ga)^2+(2-\Ga)\Ga s]\frac{c}{\delta}-(u+\frac{c}{\delta})(1-\Ga+\Ga s^{\frac{1}{2}})t^{\frac{1}{2}}\right\}.
\EQNY
So $t_u\rw 0$, as $u\rw\IF$.
\begin{align*}
\frac{\partial F_u(s,t)}{\partial s}&=\frac{\partial V_Z(s,t)}{\partial s}\cdot\frac{1}{G_u(s,t)}-\frac{V_Z(s,t)}{G_u^2(s,t)}(\frac{1}{2}\gamma u s^{-\frac{1}{2}}+\frac{c\Ga}{2\delta} s^{-\frac{1}{2}})\\
&=\frac{1}{2G_u^2(s,t) V_Z(s,t)}\left[\frac{\partial V_Z^2(s,t)}{\partial s}G_u(s,t)-V_Z^2(s,t)(\gamma u s^{-\frac{1}{2}}+\frac{c\Ga}{\delta} s^{-\frac{1}{2}})\right]\\
&=\frac{\Ga \sigma^2 s^{-\frac{1}{2}}}{4\delta G_u^2(s,t) V_Z(s,t)}\left((2-\gamma)[(1-\Ga)(u+\frac{c}{\delta})-\frac{c}{\delta}t^{\frac{1}{2}}]s^{\frac{1}{2}}-
[(1-\gamma)^2-t](u+\frac{c}{\delta})\right)
\end{align*}
so for u large enough, $F_u(\cdot,t_u)$ only can reach its maximum only at $s=1$ or $s=t_u$.
$$F_u(1,t)=\frac{\sqrt{\frac{\sigma^2}{2\delta}(1-t)}}{u+\frac{c}{\delta}(1-t^{\frac{1}{2}})},\ \ \ F_u(t,t)=\frac{\sqrt{\frac{\sigma^2}{2\delta}(1-t)}}{u+\frac{\Ga}{1-\Ga}ut^{\frac{1}{2}}+\frac{c}{\delta}(1-t^{\frac{1}{2}})}<F_u(1,t).$$
So $F_u(s,t)$ attains the maximum point only at $s=1$.
Let $\frac{\partial F_u(1,t)}{\partial t}=0$, we get $t_u=\LT(\frac{c}{\delta u+c}\RT)^2$.
\QED

\BEL\label{PST}
$\{\xi(s,t), s, t\geq 0\}$ is a zero-mean stationary Gaussian field with continuous sample paths and correlation function
$ \Ga_\xi(s,t)=\exp\LT(-(a_1|s|+a_2|t|)\RT)$ for some positive constants $a_i, i=1,2$. Further, $\lim_{u\rightarrow\IF}t_u u^2=c>0$ and $\lim_{u\rightarrow\IF}\frac{f(u)}{u}=d>0$. Then for some positive constants $S_i, i=1,2, b_i, i=1,2$,
\BQNY
&&\pk{\sup_{(s,t)\in[1-S_1u^{-2},1]\times[0,(\sqrt{S_2}u^{-1}+\sqrt{t_u})^2]}
\frac{\xi(s,t)}{1+b_1(1-s)+b_2(\sqrt{t}-\sqrt{t_u})^2}>f(u)}\\
&&\sim\mathcal{P}^{\frac{b_1}{a_1}}[0,a_1d^2S_1]\widehat{\mathcal{P}}^{\frac{b_2}{a_2},d^2a_2c}[0,a_2d^2(\sqrt{c}+\sqrt{S_2})^2]\Psi(f(u)),
\EQNY
where
$$
\widehat{\mathcal{P}}^{\frac{b_2}{a_2},d^2a_2c}[0,a_2d^2(\sqrt{c}+\sqrt{S_2})^2]:=
\E{\sup_{t\in[0,a_2d^2(\sqrt{c}+\sqrt{S_2})^2]}e^{
\sqrt{2}B(t)-t-\frac{b_2}{a_2}(\sqrt{t}-\sqrt{d^2a_2c})^2}}.
$$
\EEL
\noindent
The proof of this lemma follows along the same lines of \cite{Uniform2016}[Theorem 2.1].

\begin{figure}[H]
  \centering
  % Requires \usepackage{graphicx}
  \includegraphics[width=15cm,height=8cm]{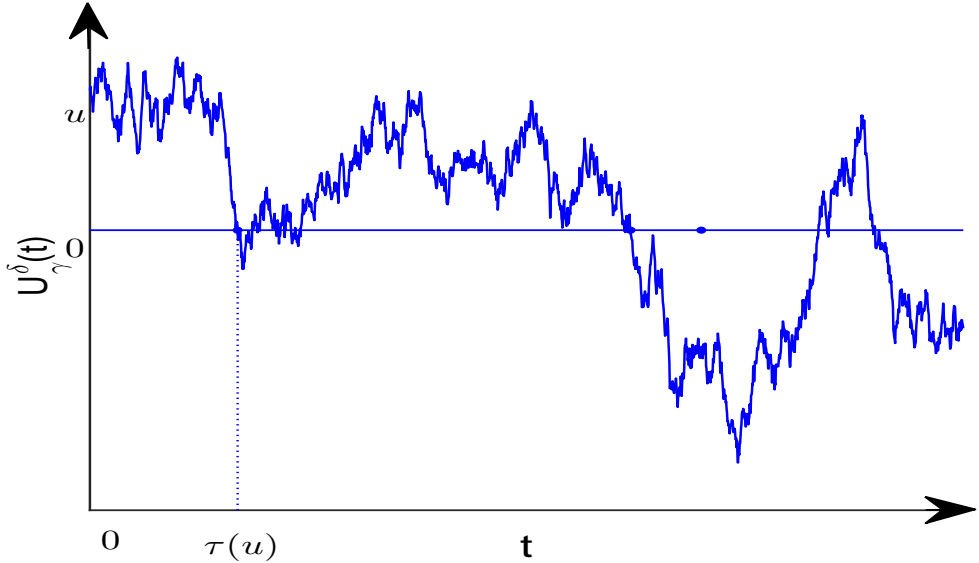}\\
  \vspace{-0.7cm}
  \caption{Ruin times}\label{fig1}
\end{figure}
Figure \ref{fig1} shows the ruin time $\tau(u)$  of a surplus process $U_\Ga^\delta(t)$.

{\bf Acknowledgement}: Thanks to Enkelejd Hashorva for his suggestions. Thanks to  Swiss National Science Foundation Grant no. 200021-166274.

\bibliographystyle{plain}

 \bibliography{gamLA}

\end{document}